\documentclass[conference]{IEEEtran}
\IEEEoverridecommandlockouts
\usepackage{cite}
\usepackage{amsmath,amssymb,amsfonts}
\usepackage{algorithmic}
\usepackage{graphicx}
\usepackage{textcomp}
\usepackage{xcolor}
\usepackage{caption}
\usepackage{subcaption}
\usepackage{url}  
\usepackage{amsmath}

\def\BibTeX{{\rm B\kern-.05em{\sc i\kern-.025em b}\kern-.08em
    T\kern-.1667em\lower.7ex\hbox{E}\kern-.125emX}}
\begin{document}

\title{
Optimizing Multi-Timestep Security-Constrained Optimal Power Flow for Large Power Grids\\
}

\author{\IEEEauthorblockN{Hussein Sharadga}
\IEEEauthorblockA{~~~~\textit{Civil, Architectural and Environmental Engineering~~~~} \\
\textit{The University of Texas at Austin}\\
Austin, Texas, USA \\
hssharadga@utexas.edu}
\and
\IEEEauthorblockN{Javad Mohammadi}
\IEEEauthorblockA{~~~~~~~~\textit{Civil, Architectural and Environmental Engineering~~~~~~~~} \\
\textit{The University of Texas at Austin}\\
Austin, Texas, USA \\
javadm@utexas.edu}
\and
\IEEEauthorblockN{Constance Crozier}
\IEEEauthorblockA{~~~~~~~~~~~\textit{ Industrial and Systems Engineering~~~~~~~~~~~~~~} \\
\textit{Georgia Institute of Technology}\\
Georgia, USA \\
ccrozier8@gatech.edu}
\and
\IEEEauthorblockN{Kyri Baker}
\IEEEauthorblockA{~~~~~~\textit{Civil, Environmental, and Architectural Engineering~~~~~~} \\
\textit{University of Colorado}\\
Boulder, Colorado, USA \\
Kyri.Baker@colorado.edu}
}
\maketitle

\begin{abstract}
This work proposes a novel method for scaling multi-timestep security-constrained optimal power flow in large power grids. The challenge arises from dealing with millions of variables and constraints, including binary variables and nonconvex, nonlinear characteristics. To navigate these complexities, techniques such as constraint relaxation, linearization, sequential optimization, and problem reformulation are employed. By leveraging these methods, complex power grid problems are solved while achieving high-quality solutions and meeting time constraints. The innovative solution approach showcases great robustness and consistently outperforms benchmark standards.
\end{abstract}

\begin{IEEEkeywords}
power grid optimization, power grid software, security-constrained optimal power flow, large-scale optimization
\end{IEEEkeywords}

\vspace{-.4cm}
\section{Introduction}

The power grid stands as one of the most significant engineering achievements of the 20th century \cite{b1}. However, optimizing power grid dispatch presents a substantial challenge in ensuring efficient and reliable energy distribution. As modern power systems grow increasingly complex, finding robust solutions becomes imperative.
To address these challenges and enhance grid resilience, the US Department of Energy (DOE) has been actively involved in initiatives focused on advancing grid technologies. These efforts aim to modernize grid operating models and software, integrating new technologies and methodologies to adapt to the complexities of the evolving grid, ensuring its reliability, resilience, and security.

The optimal power flow (OPF) problem is the core optimization challenge in grid planning and operations \cite{b2}. It involves finding the best settings for power generation, demand flexibility, energy storage, and grid control to maximize grid objectives. Security-constrained optimal power flow (SCOPF) adds grid resiliency constraints to this optimization.

Today's industry-standard SCOPF methods were developed in an era of less capable and costlier computers, when general-purpose optimization solvers were in their infancy. To simplify calculations, most often, linearizing assumptions were made, ignoring voltage and reactive power—referred to as the DC optimal power flow (DCOPF) \cite{b3}. Despite various improvements, widely-used optimization software, including production cost models, security-constrained unit commitment (SCUC), and security-constrained economic dispatch (SCED) tools, continue to rely on linear OPF assumptions similar to those in a classical DCOPF problem. There are currently no widely-adopted industry tools that use the full AC power flow equations, without linearization, while simultaneously optimizing both real and reactive power generation \cite{b4}.

Recent computational and methodological advancements suggest the potential for substantial enhancements in SCOPF methodologies. The significant leaps in computational capabilities and the evolution of optimization solvers in recent times have fuelled investigations into innovative strategies for grid operations and novel approaches to address SCOPF and related grid challenges \cite{b5}. Moving away from traditional linear approximations, current research focuses on innovative methods such as quadratic convex (QC), semi-definite (SDP), and second-order cone programming techniques \cite{b6, b7, b8, b9}. Nevertheless, it has been observed that SDP relaxation fails to yield a solution with meaningful physical implications when applied to solve OPF for numerous practical systems \cite{b10}. Additionally, the existing solvers are not as efficient as those employed for solving linear programming (LP) and second-order cone programming (SOCP) formulations, especially for larger systems and problems involving mixed integer aspects. On the other hand, QC relaxation, as in \cite{b11}, is more straightforward to implement and computationally efficient. 

Also, recent studies suggest that SCOPF problems addressing both pre-contingency and post-contingency states could benefit from emerging research in decomposition methods \cite{b13} or stochastic optimization \cite{b14} algorithms. These approaches aim to exploit a similar two-stage problem structure, presenting new avenues for improving the efficiency and effectiveness of solving complex power system optimization challenges.

In this paper, we present our innovative approach to address the challenges of security-constrained optimal power flow (SCOPF) within the framework of the US DOE's Grid Optimization (GO) initiatives. Our contributions include leveraging a two-pronged approach to provide a fast high-quality solution for the complete AC power flow equations and optimizing both real and reactive power generation simultaneously. Our novel decomposition strategy breaks down the SCOPF problem into sequential sub-problems: the DC and AC models. By employing specialized solvers for each stage and incorporating advanced optimization techniques such as quadratic convex relaxation and decomposition methods, we significantly improve the efficiency and accuracy of the solution process. 

\section{Problem Formulation}
\begin{figure*}[htbp]
\includegraphics[scale=0.7]{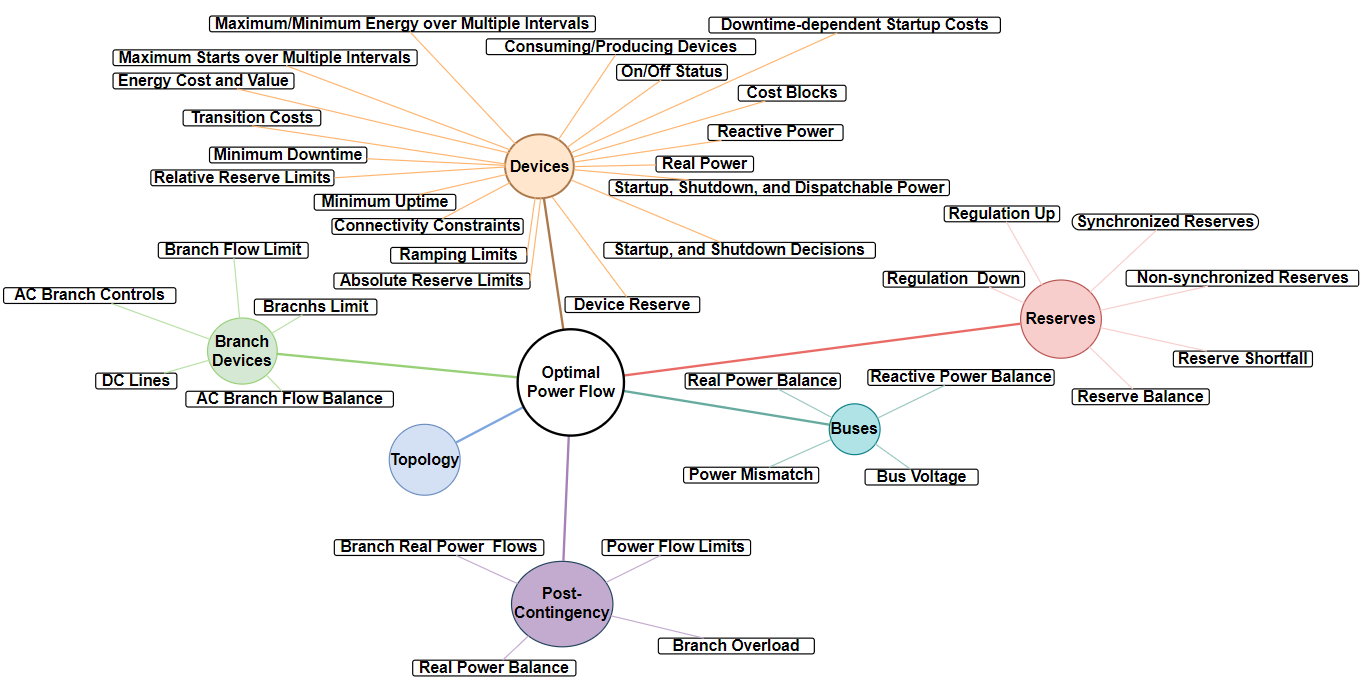}
\caption{The visual representation of Optimal Power Flow Problem formulation developed by the ARPA-E for the Grid Optimization challenge. Circles represent major components of this formulation while boxes illustrate the associated constraints.}
\label{fig1}
\end{figure*}

In this paper, we use the problem statement developed for the DOE Advanced Research Projects Agency - Energy (ARPA-E) GO challenge. 
This problem is quite extensive, with a main formulation that spans 53 pages \cite{b15}, plus an additional 22 pages of supplementary information. There are a wide range of optimization constraints and variables, including limits on real and reactive power at various bus nodes, voltage constraints, zone-specific reserves, device statuses, switching intricacies, considerations related to device downtime and uptime, the number of device starts over multiple intervals, ramping up and down constraints, and stipulations on minimum and maximum energy levels across multiple time intervals, to name a few (See Fig.~\ref{fig1}). It's important to note that all these constraints appear in both the base and contingency conditions. The contingency constraints serve to guarantee system reliability when addressing unforeseen events and outages.
The compact form of the OPF problem is:



\begin{equation}
\begin{aligned}
\max_{p_{jt},\theta_{jt}}  \, \quad& z~(p_{jt},\theta_{jt},u_{jt}, \ldots) \\
\text{s.t.} \, \quad& G~\hspace{-.13cm}(p_{jt},\theta_{jt},u_{jt}, \ldots)\hspace{-.08cm} = \hspace{-.08cm}0, \hspace{-.5cm}
\, \quad& H~\hspace{-.15cm}(p_{jt},\theta_{jt},u_{jt}, \ldots) \hspace{-.08cm}\le \hspace{-.08cm} 0
\end{aligned}
\end{equation}

Where ${z}$ is the total market surplus, $p_{jt}$ represents the power value for device ${j}$ at time ${t}$, ${\theta_{jt}}$ is the voltage angle, and $u_{jt}$ is a binary variable representing device switching decisions. In these equations, $G$ denotes the set of equality constraints, while $H$ denotes the set of inequality constraints.

\section{Problem Complexity and Testing Platform}
The network sizes span from a relatively small-scale 73 buses to a highly complex 8316 buses test system. Our dataset is comprised of 9 network models; 73-, 617-, 1576-, 2000-, 4224-, 6049-, 6708-, 6717-, and 8316-bus configurations. 
For larger-scale test cases, we need to deal with over 2 million binary variables and more than 9 million continuous variables. We will use the testing platform that is made available by the Pacific Northwest National Laboratory (PNNL). This testing platform consists of a single computing node with local memory storage, equipped with dual 32-core CPUs (64 threads). The testing criteria consists of three categories, each with unique power grid scheduling requirements:
\begin{itemize}
    \item Category 1: Optimal scheduling for an 8-hour period, with a maximum decision time of 10 minutes.
    \item Category 2: Optimal scheduling for a 48-hour window, with a maximum decision time of 2 hours.
    \item Category 3: Optimal scheduling for 7 days (168 hours), with a maximum decision time of 3 hours.
\end{itemize}

A summary of the power grid scheduling requirements for each category is given in Table \ref{tab1}.

\begin{table}[b]
\caption{Power Grid Scheduling Requirements by Category (i.e., Testing Criteria)}
\vspace{-.25cm}
\begin{center}
\begin{tabular}{|c|c|c|}
\hline
Category & Horizon & Time Constraints\\
\hline
1 & 8 hours look  ahead	& 10 mins \\
\hline
2	& 48-hour look  ahead	& 2 hours\\
\hline
3	& 7-day (168-hour) look ahead&	3 hours \\

\hline
\end{tabular}
\label{tab1}
\end{center}
\end{table}

\section{Proposed Strategy}

To address the complexities of solving large-scale SCOPF problem and reducing the solution time, we propose breaking the original problem into two sequential sub-problems: the DC and AC modules.
The initial subproblem involves the DC assumption, wherein binary and continuous are optimization variables. The binary variables returned by the first subproblem are fixed, while the continuous variables are discarded. 

In the AC module, our aim is to include all relevant constraints to accurately represent a realistic grid. However, in the DC module, we choose not to include all constraints to reduce the computation time of the first subproblem (module). Nevertheless, we do incorporate constraints essential for ensuring the feasibility of the AC subproblem (module), e.g., we include ramping limits. Without these limits in the DC module, it might prompt abrupt device shutdowns. Consequently, the AC module, which accounts for ramping limits, may identify these shutdowns as infeasible to ramp down.

\begin{figure*}
\centerline{\includegraphics[scale=0.17]{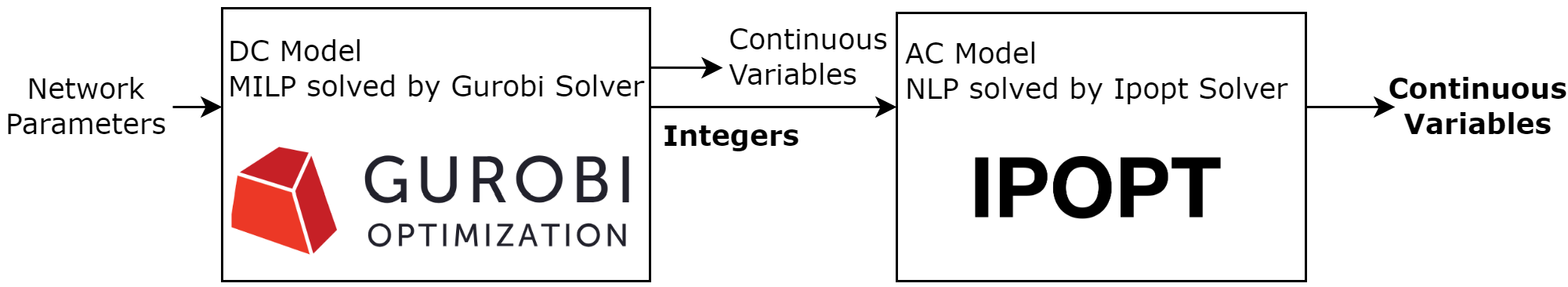}}
\caption{Our two-pronged solution approach which includes a DC and an AC module.}
\vspace{-.3cm}
\label{fig2}
\end{figure*}

%
The DC model is formulated as a mixed-integer linear problem (MILP) and is solved using the Gurobi solver. On the other hand, the AC model is a nonlinear (NLP) and nonconvex problem, which is efficiently solved by the IPOPT solver. 
By decomposing the problem and leveraging specialized solvers for the sub-problem of each module, we increase computational effectiveness and overall efficiency. Fig.~\ref{fig2} illustrates our proposed two-stage solution approach and it what follows we will discuss the novelties of each module.


\vspace{-.09cm}
\subsection{Managing Computational Workload in the DC Model}

In managing the computational workload within the DC model, we implement a series of techniques to reduce complexity and improve efficiency.

\paragraph{Eliminating Quadratic Terms by Reformulation} While quadratic terms can typically be linearized or relaxed by adding more variables and constraints, we achieved the transformation of these terms into linear expressions without introducing new variables or constraints. This artful reformulation significantly reduces computation burden. For example, the energy balance constraint is represented as \cite{b4}:

\small
\begin{equation}
\sum_{j \in j_i^{cs,pr}} p_{jt} u_{jt}^{on}+
\sum_{j \in j_i^{sh}} p_{jt} +
\sum_{j \in j_i^{fr}} p_{jt}^{fr}+
\sum_{j \in j_i^{to}} p_{jt}^{to}=p_{it}\label{eq7}
\end{equation}
\normalsize

Here, ${p_{jt}}$ in the first term of \eqref{eq7} denotes the device power, as a generator (denoted by the superscript 'pr') or a consumer (denoted by the superscript 'cs'). The power limit is given as:

\small
\begin{equation}
p_{jt, min} \le p_{jt} \le p_{jt, max} \quad \forall j \in j_{i} ^{cs,pr}\label{eq8}
\end{equation}
\normalsize

The  ${u_{jt}^{on}}$ is the device on/off status. The quadratic term ${p_{jt} u_{jt}^{on}}$ can be linearized by introducing a new variable and four constraints. However, our approach reformulates the problem into a linear format without the need for additional variables or constraints. Equations \eqref{eq7} and \eqref{eq8} are updated as:

\small
\begin{equation}
\sum_{j \in j_i^{cs,pr}} p_{jt}+
\sum_{j \in j_i^{sh}} p_{jt} +
\sum_{j \in j_i^{fr}} p_{jt}^{fr}+
\sum_{j \in j_i^{to}} p_{jt}^{to}=p_{it}\label{eq9}
\end{equation}
\normalsize

\begin{equation}
p_{jt, min} u_{jt}^{on}\le p_{jt} \le p_{jt, max} u_{jt}^{on} \quad \forall j \in j_{i} ^{cs,pr}\label{eq10}
\end{equation}

\paragraph{Nonconvex Constraints Relaxation} This involves the utilization of the substitution method and the relaxation of the equality constraint within the nonlinear constraints into an inequality constraint. For instance, shunt power (denoted by the superscript 'sh') is represented as \cite{b4}:



\small
\begin{equation}
p_{jt}=g_{jt}^{sh} v_{it}^{2}\label{eq2} \quad \forall j \in j_{i} ^{sh}
\end{equation}
\normalsize

Equation \eqref{eq2} is non-linear. In an effort to eliminate the non-linearity, we will first substitute \eqref{eq2} into \eqref{eq9}:

\small
\begin{equation}
\sum_{j \in j_i^{cs,pr}} p_{jt}+
\sum_{j \in j_i^{sh}} g_{jt}^{sh} v_{it}^{2} +
\sum_{j \in j_i^{fr}} p_{jt}^{fr}+
\sum_{j \in j_i^{to}} p_{jt}^{to}=p_{it}\label{eq3}
\end{equation}
\normalsize

However, the resulting constraint is nonconvex. 
The new constraint can be transformed into a convex constraint by relaxing the equality constraint into an inequality constraint:

\small
\begin{equation}
\sum_{j \in j_i^{cs,pr}} p_{jt}+
\sum_{j \in j_i^{sh}} g_{jt}^{sh} v_{it}^{2} +
\sum_{j \in j_i^{fr}} p_{jt}^{fr}+
\sum_{j \in j_i^{to}} p_{jt}^{to}\le p_{it}\label{eq4}
\end{equation}
\normalsize

Because the right side of \eqref{eq4} represents the power mismatch, $p_{it}$, and due to a significant penalty imposed on this mismatch, the solver strives to minimize it. The mismatch will be minimized when it equals the value on the left side of \eqref{eq4}.

\paragraph{Convex Constraints Linearization}  First, we reformulate the convex constraints into a quadratic problem. Next, we relax the quadratic constraints into linear problems using McCormick envelope relaxation \cite{b7}. This reduction in complexity leads to a significant 
reduction in computation time for small networks and even greater time savings for larger networks. Specifically, the branch flow limit and penalty are represented as a second-order cone problem 
taking the following format:


\small
\begin{equation}
||A x_i+B y_i||_2 \le Cz_i+D\label{eq5}
\end{equation}
\normalsize

Here, the left side is the norm function which is a convex function. Since the right side is linear, \eqref{eq5} is convex. We reformulated \eqref{eq5} into quadratic by squaring both parts as,
\small
\begin{equation}
A^2 x_i^2 + B^2 y_i^2 + 2AB x_i y_i  \le C^2 z_i^2 +D^2+2 C D z_i\label{eq6}
\end{equation}
\normalsize
While \eqref{eq6} is nonconvex, the bilinear and quadratic terms can be relaxed into linear using a McCormick envelope. 

\paragraph{Sparse Constraint-matrix} The MILP problem, given this formulation, is always feasible. Nevertheless, introducing slack variables aids in faster convergence for the solver. These slack variables appear in only a few constraints, resulting in a sparser constraint matrix with more zeros, which facilitates quicker convergence for the solver.

\paragraph{Reduction of Lengthy Linear Expressions} A significant number of devices are connected to the same bus, resulting in an extensive power balance constraint in the form of a lengthy linear expression. This complexity can lead to numerical issues for the solver. To mitigate this, we introduce new variables. Even though this increases optimization variables, our adjustment enhances convergence speed.

\paragraph{Reformulating Downtime-Dependent Startup Costs} We formulated downtime-dependent startup costs using seven different models, which include three conventional models, the polyhedron representation transformation algorithm, and three versions of the clique constraint \cite{b16}. 
Rather than selecting the model with the fewest variables and constraints, our choice is based on compatibility with the Gurobi workflow. Specifically, we opted for the clique constraint model due to its alignment with Gurobi's pre-solver, which uses problem structure for eliminating variables and constraints.

\paragraph	{Max/Min Constraint Functions Reformulation} The maximum and minimum equality constraints can be reformulated in two ways: (1) as a linear model or (2) as an MIP problem. Our analysis suggests that the MIP formulation aligns more effectively with Gurobi's structure and workflow.

\paragraph{Heuristic Line Switching Algorithm} In our problem, the decision regarding line switching is based on binary optimization variables. The introduction of additional binary variables can significantly increase the complexity of the problem-solving process. When line switches are treated as variables rather than being fixed at their initial conditions, it transforms the problem from a linear one into a quadratic one. Linearizing the problem, introduces additional variables and constraints, which, in turn increase the computational burden. To mitigate this challenge, we have developed a heuristic line-switching algorithm. This algorithm effectively disconnects branches in parallel with reactance or impedance values differing from those in their immediate vicinity.

\subsection{Managing Computational Workload in the AC Model}

This subsection's efforts are focused on managing the computational workload in the AC module using various techniques. Notably, we handle the entire AC model without any linearization or relaxation. This ensures a comprehensive and accurate representation of the power grid system. The details of these techniques are presented below:

\paragraph{Handling Controllable and Uncontrollable Load} Every device has a set of cost blocks ${(m\in M)}$ [4]. We investigated every device to identify pseudo-variables that cannot be controlled. This approach reduces the number of variables associated with the cost blocks and eliminates the need for introducing binary variables. More precisely, the power is split into bid blocks ${(m\in M)}$, and a cost is applied to each block. Therefore, the total power of the device is the sum of the power values in each block:

\small
\begin{equation}
p_{jt} = \sum_{m \in M_{jt}}  p_{jtm} \quad  \forall j \in j_{i} ^{cs,pr}\label{eq11}
\end{equation}
\normalsize

For instance, if the minimum power value of the device is higher than the power in the first block ${(m=1)}$, then the power value in the first block is uncontrollable. Hence, the power value in the first block is equal to its maximum power.

\paragraph{Post-processing Reserves Calculations} Reserves are not incorporated into the main optimization framework. Instead, after IPOPT completes the simulation, power values are processed to meet the reserve requirements. This modification significantly reduces the problem's complexity. Another approach we considered was incorporating the reserve constraints in the DC section and then using the obtained power values as the maximum power for the AC section.

\paragraph{Exact Hessian and Jacobian Matrix} We compute the first and second derivatives of the Hessian and Jacobian analytically. Without these matrices, the solver resorts to less efficient numerical methods for estimation. Supplying the matrices analytically expedites IPOPT's convergence but necessitates a more hands-on code development, rather than relying on IPOPT's numerical techniques.

\paragraph{Providing Hessian and Jacobian Structure} When supplying the exact Hessian and Jacobian matrices, we provide information exclusively for non-zero elements along with their respective locations, as indicated by the structure. This approach eliminates the need to provide information for zero-valued elements, effectively reducing the computational workload and the volume of stored data. As a result, it contributes to faster convergence.

\paragraph{Utilizing Vectorized Form and Sparse Matrix in Coordinate format} This minimizes the required storage capacity. However, when not using the sparse matrix in coordinate format, we ran out of memory.

\begin{figure}[b]
\vspace{-.5cm}
\centerline{\includegraphics[scale=0.15]{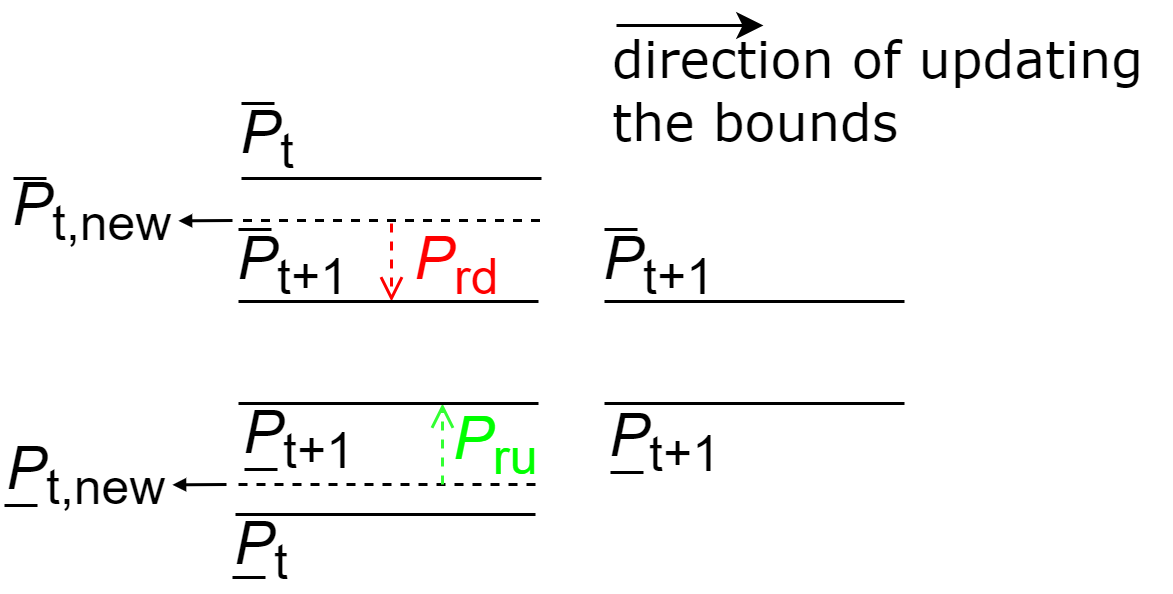}}
\caption{Bound update in stage 1 to prevent ramping violations in sequential solving.}
\label{fig3}
\end{figure}

\begin{figure}[t]
\vspace{-.3cm}
\centerline{\includegraphics[scale=0.15]{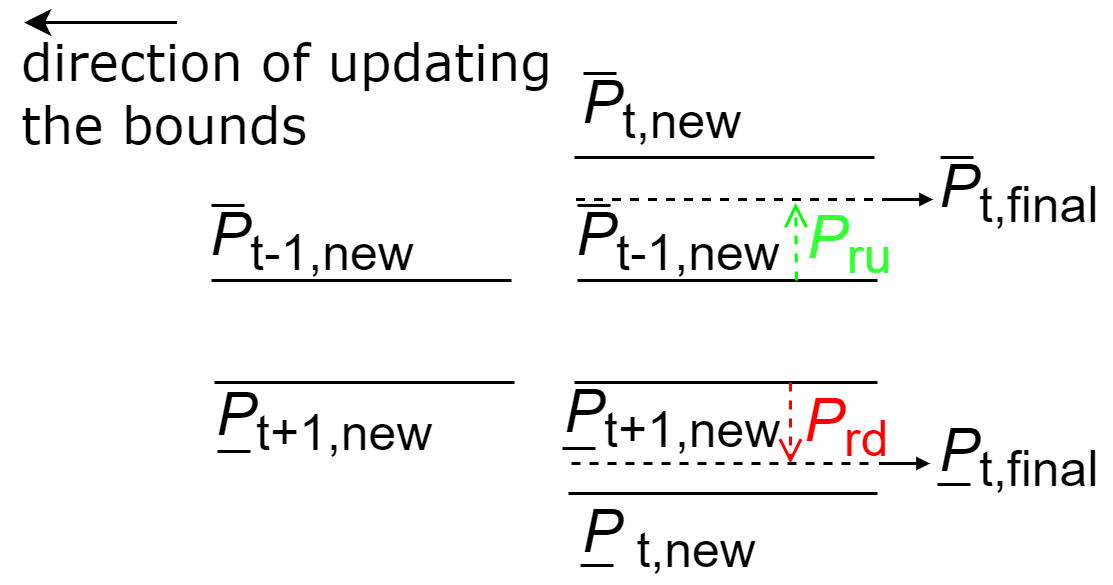}}
\caption{Update bounds in 2nd stage to avoid ramping violations during sequential solving.}
\label{fig4}
\end{figure}

\paragraph{Receding Ramping Bounds} This approach facilitates solving the problem through sequential optimization. The problem is addressed sequentially for large networks because it is not feasible to solve these extensive networks within the given time constraints if we attempt to tackle the entire problem without breaking it down into steps. Therefore, we recursively update the bounds to ensure that we stay within the ramping limits when solving each step individually. The process of updating the bounds is undertaken in two stages, as illustrated in Figs.~\ref{fig3} and \ref{fig4}.
The bounds are updated in stage 1 using the following equations as we move forward:
\small
\begin{align}
\overline P_{t, new} = \min [ \overline P_{t}, \overline P_{t+1} +P_{rd} ]\label{eq12}\\
\underline P_{t, new} = \max [ \underline P_{t}, \underline P_{t+1} - P_{ru} ]\label{eq13}
\end{align}
\normalsize
The second stage bounds are updated as we move backward:
\small
\begin{align}
\overline P_{t, final} = \min [ \overline P_{t, new}, \overline P_{t-1, new} +P_{ru} ]\label{eq14}\\ 
\underline P_{t, final} = \max [ \underline P_{t, new}, \underline P_{t-1, new} - P_{rd} ]\label{eq15}
\end{align}
\normalsize

\subsection{Efficient Resource Sharing between DC and AC Modules}
Addressing the complex nonconvex nature of the AC model demands a considerable amount of computational resources compared to the MILP problem. To optimize our two-pronged approach, we divide our time allocation as follows: One-third of  the computational resources are dedicated to solving the DC model, while the remaining two-thirds are devoted to tackling the challenges posed by the AC model. This balanced allocation ensures efficient coordination between the DC and AC modules, allowing us to solve these problems in time.
\subsection{Sequential Optimization}
Solving the problem in a single shot is impractical given the problem scale, typically those with more than 2000 buses. In such cases, the control horizon is subdivided into stages, where the outcomes of one stage serve as the initial conditions for the subsequent one.
In the context of the mixed-integer programming (MIP) problem, devices are permitted to undergo two switching events: the first occurs at the initiation of the control horizon. However, a crucial constraint to consider is the minimum uptime requirement, which necessitates that the device remains operational for a certain duration before it can be turned off again. As a result, the second switching event can only occur once this minimum uptime has been satisfied.

To ensure compliance with ramping limits and prevent violations, the bounds are adjusted in a recursive manner before executing the sequential optimization process. 
This method enables the efficient handling of intricate scenarios in large-scale networks while adhering to critical operational constraints.
These tactics collectively contribute to a streamlined and computationally efficient solution process. By integrating these techniques, our approach consistently continues to excel in providing an effective solution to power grid optimization.

\section{Performance Evaluation}

\begin{figure}
\vspace{-.2cm}
\centerline{\includegraphics[scale=0.7]{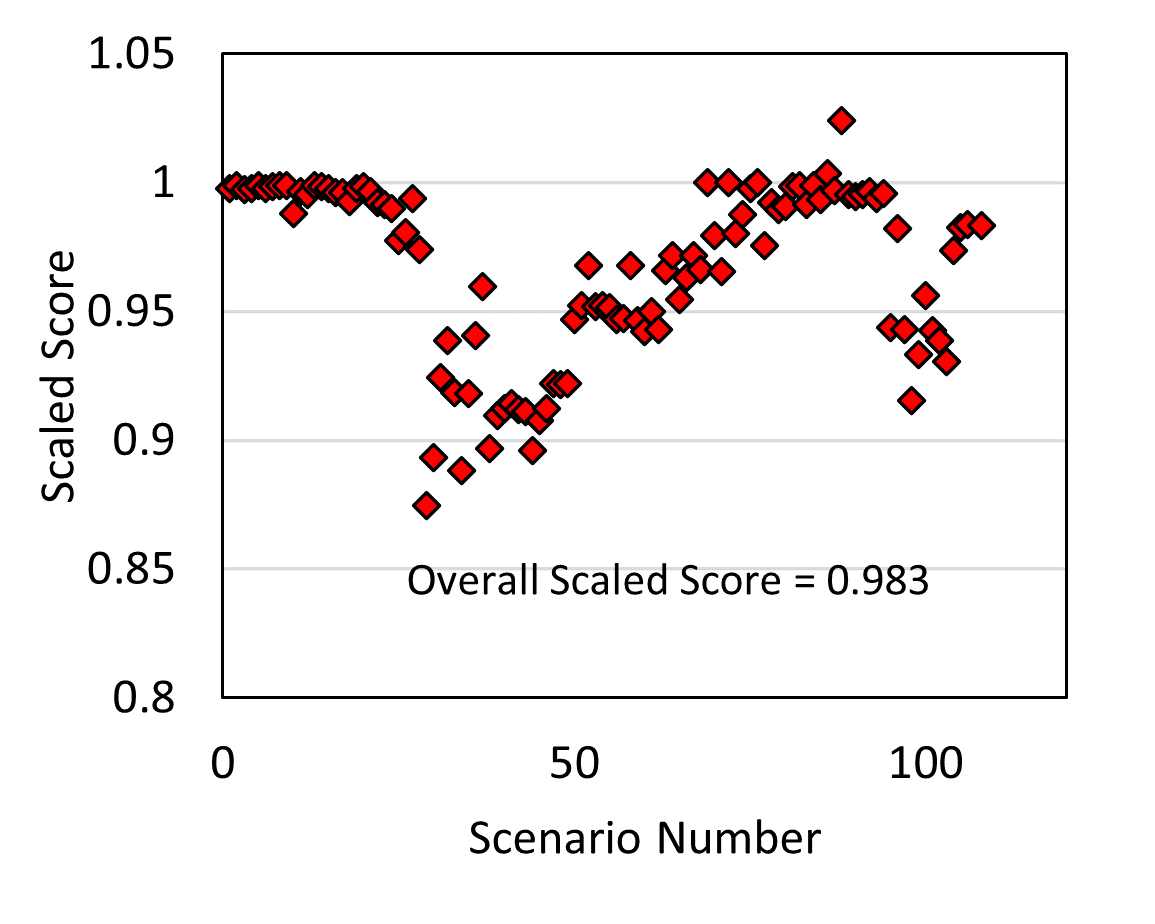}}
\vspace{-.35cm}
\caption{Performance of the proposed solution strategy in scenarios set 3, evaluated against the best available solution. The scaled score represents the ratio to the best score.}
\label{fig5}
\vspace{-.3cm}
\end{figure}

\begin{figure} 
     \centering
     \begin{subfigure}[b]{0.48\textwidth}
         \centering
         
         \vspace{-.2cm}
         \includegraphics[width=0.7\textwidth]{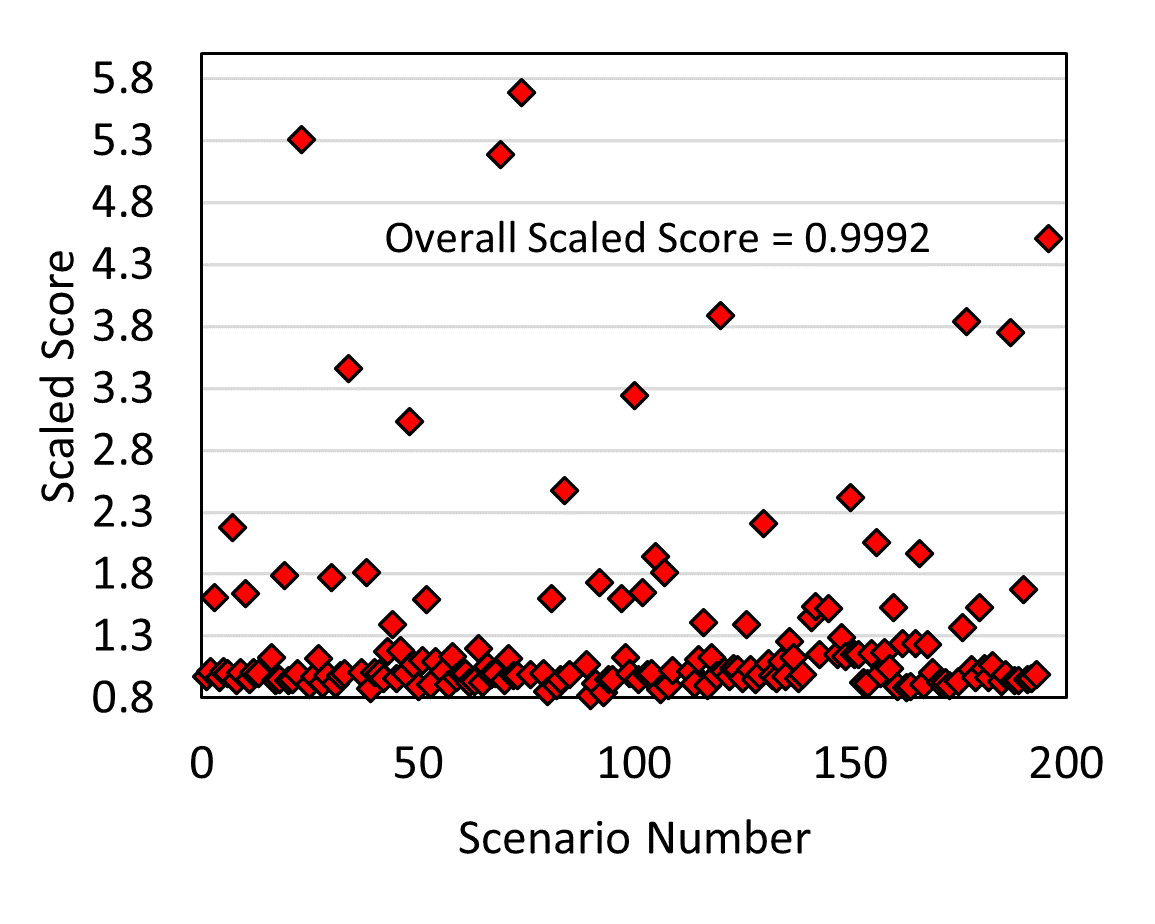}
         \vspace{-.47cm}
         \caption{}
         \label{fig6.a}
     \end{subfigure}
     \hfill
     \vspace{-.09cm}
     \begin{subfigure}[b]{0.48\textwidth}
         \centering
         \includegraphics[width=0.7\textwidth]{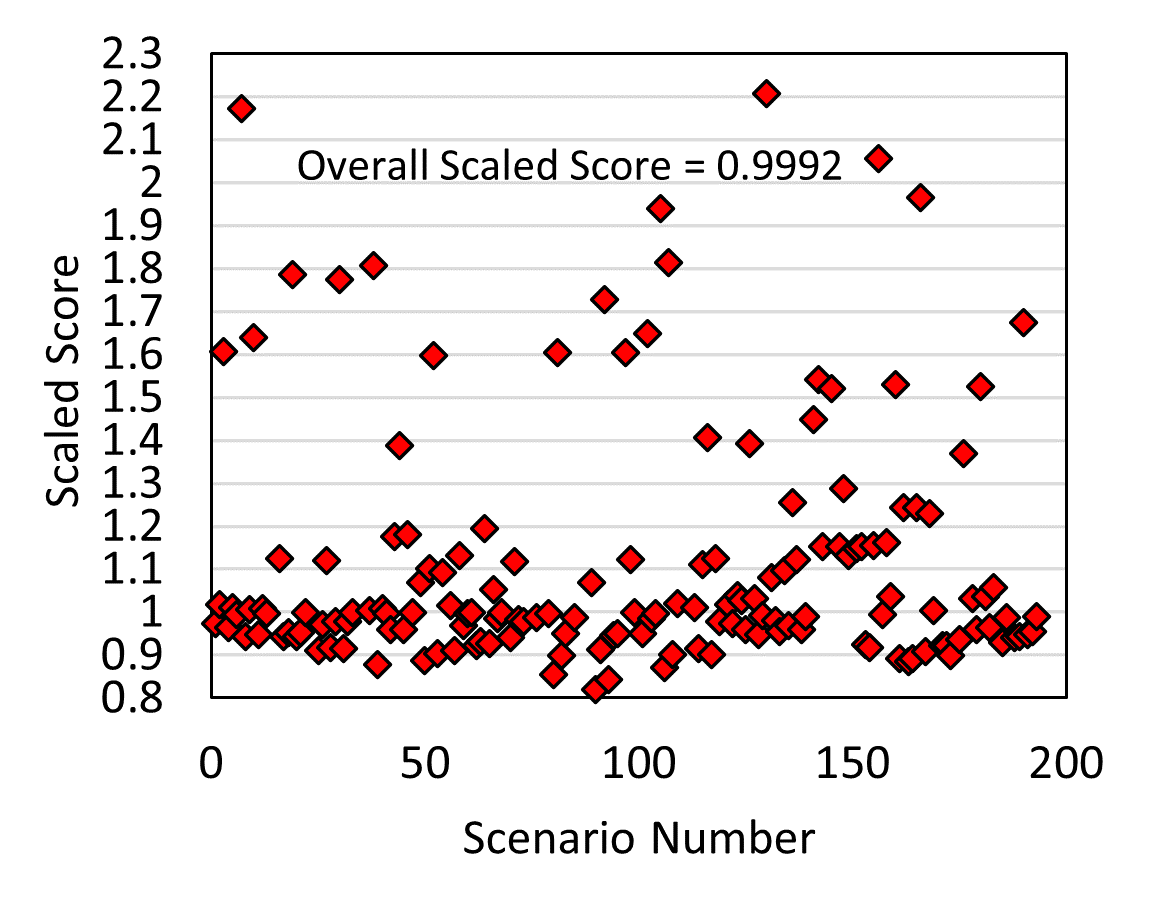}
         \vspace{-.47cm}
         \caption{}
         \label{fig6.b}
     \end{subfigure}
     \vspace{-.15cm}
        \caption{Performance of the proposed solution strategy in (a) scenarios set 2 and (b) zoomed-in view of scenarios set 2, evaluated against the best available solution.}
        \label{fig6}
\end{figure}

Our method was tested using over 300 scenarios. The code was evaluated based on two criteria: the total score and the number of best scores achieved.	The scaled scores for scenarios set 3 and scenarios set 2 are displayed in Figs.~\ref{fig5} and \ref{fig6}. The scaled score is calculated as the achieved score divided by the highest score recorded among available codes. A scaled score of 1.0 indicates the highest possible achievement, while a score of 0.99 signifies a ${1\%}$ deviation from the top performance.

After reviewing Fig.~\ref{fig5} and Fig.~\ref{fig6}, it's evident that the scaled score for the majority of scenarios is higher than 0.95, with an overall scaled score exceeding 0.98. This demonstrates the robustness of the proposed strategy in effectively handling various scenarios.
More precisely, in scenarios set 3, the code was tested in 137 scenarios, while in scenarios set 2, it was tested in 196 scenarios. In scenarios set 3, the scaled score was at least 0.99 in ${42\%}$ of the scenarios, whereas in scenarios set 2, this was the case in ${53\%}$ of the scenarios. Importantly, the scaled score exceeded 0.95 in ${70\%}$ of the scenarios for both scenarios set 3 and scenarios set 2.  The overall scaled scores for both scenarios set 3 and scenarios set 2 are 0.983 and 0.9992, respectively. A detailed distribution of these score ranges can be found in Table \ref{tab2}.



\begin{table}[b]
\caption{Distribution of Scaled Scores for Set 3 and 2}
\vspace{-.2cm}
\begin{center}
\begin{tabular}{|c|c|c|}
\hline
\textbf{Range} &\multicolumn{2}{|c|} {\textbf{${\%}$ Scenarios}}\\
\cline{2-3}
 & \textbf{Set 3} & \textbf{Set 2}\\
\hline
Scaled Score ${\ge}$ 0.99 	& 42 & 53\\
\hline
${0.99>}$ Scaled Score ${\ge}$ 0.97 & 16 & 10\\
\hline
${0.97>}$ Scaled Score ${\ge}$ 0.95 & 12 & 9 \\
\hline
${0.95>}$ Scaled Score ${\ge}$ 0.93 & 12 &9 \\
\hline
${0.93>}$ Scaled Score ${\ge}$ 0.90 & 13 & 11 \\
\hline
 Scaled Score ${<}$ 0.9 	& 5 & 8\\
\hline
\end{tabular}
\label{tab2}
\end{center}
\end{table}

\begin{figure}
\centerline{\includegraphics[scale=0.7]{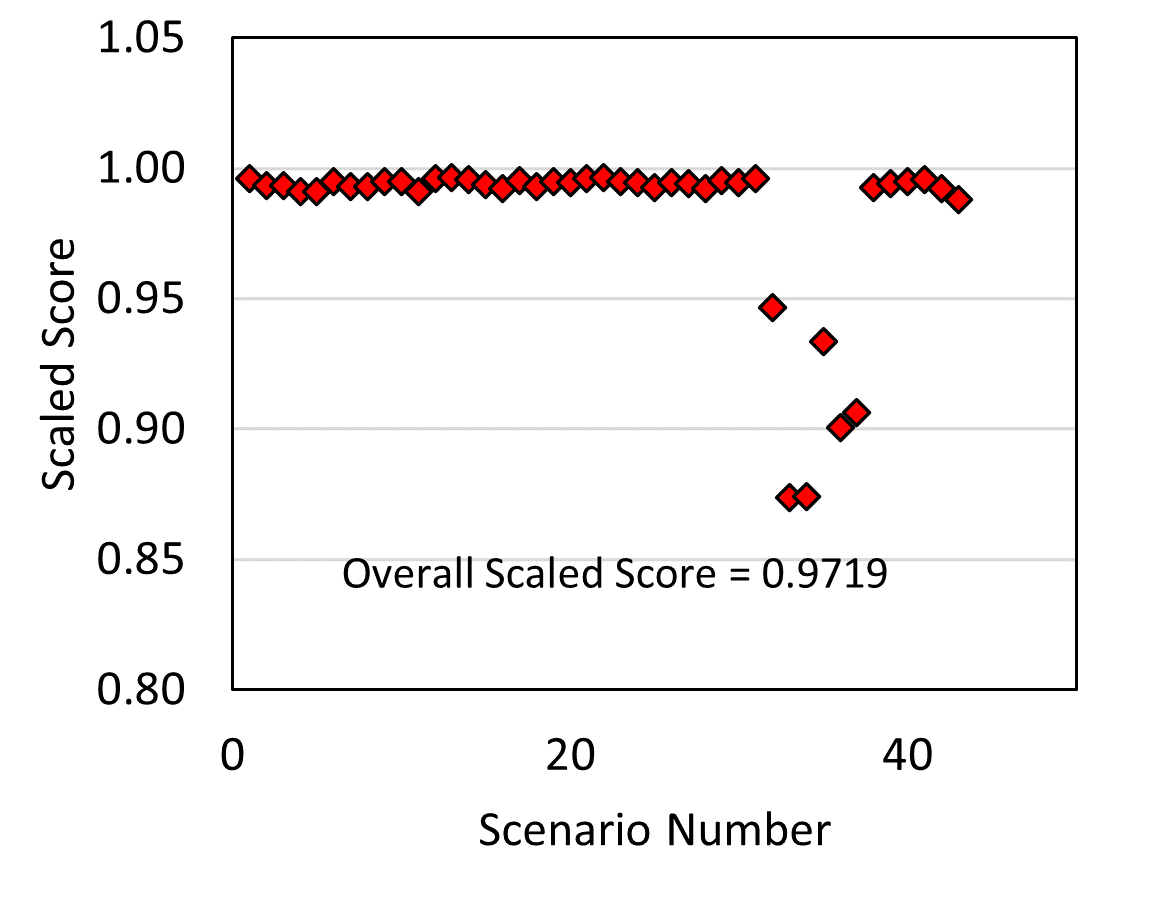}}
\vspace{-.5cm}
\caption{Performance of the proposed  strategy on 6708-bus industrial network, evaluated against the best available solution.}
\label{fig7}
\vspace{-.5cm}
\end{figure}

The data set encompassed a variety of network models, each differing in size and complexity. These models included configurations with 73, 617, 1576, 2000, 4224, 6049, 6708, 6717, and 8316 buses, showcasing a wide range of network scales. The 6708-bus network is the only industrial network; the rest are synthesized networks. In the industrial network scenarios, totaling approximately 45 (as shown in Fig.~\ref{fig7}), the scaled score was nearly optimal in about 39 scenarios. Among the remaining 6 scenarios, 4 exhibited scaled scores between 0.95 and 0.9, while the other 2 scenarios fell just below the 0.9 threshold. This demonstrates our solution's performance in solving this large-scale industrial real network.



\section{Conclusion}
\vspace{-.1cm}
In conclusion, addressing the intricate challenge of optimizing the power grid requires a highly effective strategy. The endeavor to scale multi-timestep security-constrained optimal power flow in large power grids demands innovative solutions due to the involvement of millions of variables and constraints, including binary variables and nonconvex, nonlinear characteristics. Our solution approach, grounded in techniques such as constraint relaxation, linearization, sequential optimization, and problem reformulation, has proven to be efficient in navigating the complexities of the problem while adhering to time constraints and maintaining optimality. The strategic division of the problem into sequential sub-problems, namely the DC model and the AC model, along with the careful management of computational workload, proved essential in enhancing the efficiency and effectiveness of our approach. The performance evaluation highlights the robustness of our strategy in addressing the intricate challenges inherent in power grid optimization on a large scale.

\section*{Acknowledgment}
\vspace{-.05cm}
This work is financially supported by the US ARPA-E (\#DE-AR0001646). We used ChatGPT to increase readability.


\end{document}